\documentclass[a4j,12pt]{article} 
\usepackage{amsmath,amssymb}

\newcommand{\Section}[1]{%
\renewcommand{\thesection}{\S\arabic{section}}
\section{#1}
\renewcommand{\thesection}{\arabic{section}}
\setcounter{equation}{0}}
\newcommand{\qed}{\hfill \hbox{\rule[-2pt]{4pt}{7pt}}}
\newcommand{\proof}{{\hspace*{0.4cm} {\it Proof}.\ \enskip}}
\newtheorem{lem}{Lemma}[section]
\newtheorem{thm}[lem]{Theorem}  %あるかも
\newtheorem{prop}[lem]{Proposition}

\newtheorem{rem}[lem]{Remark}

\newcommand{\tr}{{\rm tr}}    

\newcommand{\Q}{{\Bbb Q}}
\newcommand{\Z}{{\Bbb Z}}

\newcommand{\C}{{\Bbb C}}

\newcommand{\N}{{\Bbb N}}

\newcommand{\set}[2]{\left\{\left.#1\vphantom{#2}\:\right\vert\:#2\right\}}
\newcommand{\wt}{\widetilde}

\newcommand{\lam}{{\lambda}}
\newcommand{\Lam}{{\Lambda}}

\newcommand{\alp}{{\alpha}}  %2003
\newcommand{\slit}{\vspace{5mm}}
\newcommand{\mslit}{\vspace{3mm}}

\newcommand{{\any}}{{}^\forall}
\newcommand{{\is}}{{}^\exists}
\newcommand{{\st}}{\; {\rm s.t.}\; }

\newcommand{\Llra}{{\Longleftrightarrow}}

\newcommand{\abs}[1]{\left\vert{#1}\right\vert}  %in mymacro.tex

\newcommand{\nvert}{\mbox{$\not\vert$}}

\newcommand{\mapright}[1]{\displaystyle{
   \smash{\mathop{\hbox to 1cm{\rightarrowfill}}^{#1}}}}

\newcommand{\gyaddots}{%
\setlength{\unitlength}{1mm}
\begin{picture}(5,3.5)(-2,-1.5)
\put(0,0){$\cdot$}
\put(2,1.5){$\cdot$}
\put(-2,-1.5){$\cdot$}
\end{picture}}
\newcommand{\gyakuddots}{\smash{\lower0.3ex\hbox{\gyaddots}}}

\renewcommand{\slit}{{\vspace{5mm}\noindent}}
\renewcommand{\mslit}{{\vspace{3mm}\noindent}}
\newcommand{\mmslit}{{\vspace{1mm}\noindent}}

\newcommand{\gen}[1]{\langle #1 \rangle}
\newcommand{\hdeg}{{\rm hdeg}}
\newcommand{\hterm}{{\rm hterm}}
\newcommand{\htermcoef}{{\rm htermcoeff}}
\newcommand{\ord}{{\rm ord}}

\allowdisplaybreaks[1]   %%% display の途中の改ページを許す．許さない行は \\* とせよ

\begin{document}
\title{Zeta functions of finite groups \\by enumerating subgroups}
\author{Yumiko Hironaka\\
{\small Department of Mathematics, }\\
{\small Faculty of Education and Integrated Sciences, Waseda University}\\
{\small Nishi-Waseda, Tokyo, 169-8050, JAPAN}}

\date{}
\maketitle

%\hfill \today

\begin{center}\begin{minipage}{12cm}{\small 
\hspace*{5cm}abstract

For a finite group $G$, we consider the zeta function $\zeta_G(s) = \sum_{H} \abs{H}^{-s}$, where $H$ runs over the subgroups of $G$. First we give simple examples of abelian $p$-group $G$ and non-abelian $p$-group $G'$ of order $p^m, \; m \geq 3$ for odd  $p$ (resp. $2^m, \; m \geq 4$) for which $\zeta_G(s) = \zeta_{G'}(s)$. Hence we see there are many non-abelian groups whose zeta functions have symmetry and Euler product, like the case of abelian groups. On the other hand, we show that $\zeta_G(s)$ determines the isomorphism class of $G$ within abelian groups, by estimating the number of subgroups of abelian $p$-groups. Finally we study the problem which abelian $p$-group is associated with a non-abelian group having the same zeta function.
} %small
\end{minipage}
\end{center} %\end{abstract}

\renewcommand{\thefootnote}{{}}

\footnotetext{

This research is partially supported by Grant-in-Aid for Scientific Research (C): 24540031.}

\setcounter{section}{-1}

%   
%\Section{Introduction}
%\input{intro}
%

%\input{fGZetaS1}
\Section{}%introduction}
Let us consider zeta functions associated to finite groups by enumerating their subgroups. 
For a group $G$, one may consider the following zeta functions if they convergent in some region
\begin{eqnarray}
\zeta_G(s) = \sum_{n=1}^\infty \frac{a_n(G)}{n^s}, \quad
\zeta_G^*(s) = \sum_{n=1}^\infty \frac{a_n^*(G)}{n^s},
\end{eqnarray}
where $s \in \C$ and $a_n(G)$ (resp. $a_n^*(G)$) is the number of subgroups of order $n$ (resp. index $n$ ) in $G$. 
If $G = \Z$, $\zeta_G^*(s)$ is nothing but the Riemann zeta function, and this is an example of Solomon's zeta functions for lattices. It seems to be natural to consider $\zeta_G^*(s)$ in general and there are many studies in that direction. As for finitely generated (infinite) groups, there are interesting results by F.~J.~Grunewald, M.~P.~F.~du Sautoy, et al. (e.g. \cite{GSS}, \cite{Sau}, \cite{SG}).  

\mmslit
In this paper we concentrate on finite groups and study $\zeta_G(s)$. There is no essential difference between the two zeta functions, since there is a simple relation 
\begin{eqnarray}  \label{zeta-zeta} 
\zeta_G(s) = \frac{1}{\abs{G}^s}\, \zeta_G^*(-s), \quad (\mbox{for a finite group $G$}).
\end{eqnarray}
Moreover, the duality for finite abelian groups yields that $\zeta_G(s) = \zeta_G^*(s)$, hence one has 
\begin{eqnarray}  \label{sym of zeta} 
\zeta_G(s) = \frac{1}{\abs{G}^s}\, \zeta_G(-s), \quad(\mbox{for a finite abelian group $G$}).
\end{eqnarray}
   
\mmslit
In the following we consider only finite groups. It is natural to ask

\mslit
\qquad [A] whether $\zeta_G(s)$ determines the isomorphism class of $G$ or not, 

\mslit
\qquad [B] whether the symmetry (\ref{sym of zeta}) is essential for abelian groups or not. 

\mslit
In \S 1, we will construct examples of abelian groups and non-abelian groups that have the same zeta functions, which shows Questions [A] and [B] are negative in general. On the other hand, we will show Question [A] is affirmative within abelian groups in \S 2.    
For abelian groups, Question [A] is reduced to enumerating subgroups of abelian $p$-groups and combinatorial studies, and there is a close relation to automorphic forms and local densities of matrices (cf. Remark~\ref{rem2}).     
Finally in \S 3, we give a partial answer for the problem which abelian $p$-group is associated with a non-abelian group having the same zeta function.

\vspace{2cm}
\Section{}
The zeta function $\zeta_G(s)$ is a polynomial with non-negative integral coefficients of $p^{-s}$ with various $p$, where and henceforth we denote by $p$ a prime number.
We start with an easy observation.

\begin{prop}
{\rm (1)} The zeta function $\zeta_G(s)$ determines the order $\abs{G}$ of $G$.

{\rm (2)} If $\abs{G}$ and $\abs{G'}$ are coprime, then $\zeta_{G\times G'}(s) = \zeta_G(s) \times \zeta_{G'}(s)$.

{\rm (3)} The fact only $1$ appears as coefficients of $\zeta_G(s)$ is equivalent to that $G$ is cyclic.
\end{prop}

\proof
The statements (1) and (3) are clear. Under the assumption of (2), any subgroup of $G \times G'$ is written uniquely as a direct product of the form $H \times H'$ for some subgroups $H$ of $G$ and $H'$ of $G'$. \qed

\slit
We say $\zeta_G(s)$ has an Euler product if it can be decomposed as
\begin{eqnarray} \label{Euler prod}
\zeta_G(s) = \prod_{p \mid \abs{G}}\, f_p(p^{-s}), \quad f_p(X) \in \C[X].
\end{eqnarray}
Then we have

\begin{prop} \label{prop Euler prod}
$\zeta_G(s)$ has an Euler product if and only if $G$ is a nilpotent group, i.e. $G$ is a direct product of its Sylow $p$-subgroups.
The Euler product (\ref{Euler prod}) is uniquely determined under the condition $f_p(0) = 1$ for any $p$ dividing $\abs{G}$.
\end{prop}

\proof
Assume $\zeta_G(s)$ has an Euler product of the shape (\ref{Euler prod}).
Since $\prod_{p}f_p(0) = 1 (= a_1(G))$, we may assume the constant term of $f_p(X)$ is equal to $1$ for any prime divisor $p$ of $\abs{G}$.
For such a $p$, the coefficient of $X^m$ in $f_p(X)$ must coincide with $a_{p^m}(G)$, hence it is a non-negative integer.  Denote by $c_p$ the leading coefficient of $f_p(X)$. Then $\prod_{p\mid \abs{G}}\, c_p = 1$ yields $c_p = 1$ for any $p$, which means there is only one normal Sylow $p$-subgrop in $G$ for each $p$ and $G$ is nilpotent. The converse and the uniqueness of the Euler product are clear.
\qed

\slit
By the above proof, we see that the decomposition (\ref{Euler prod}) is unique if one takes the constant terms of the Euler factors as $1$ and actually given by polynomials with non-negative integral coefficients.

\mslit
Here we give an example which does not have an Euler product.     
Take an odd integer $n (\geq 3)$ and consider the dihedral group of order $2n$
\begin{eqnarray*}
D_n = \gen{\sigma, \tau \mid \sigma^n = \tau^2 = 1, \; \tau\sigma\tau = \sigma^{-1}}.
\end{eqnarray*}
Then
\begin{eqnarray*}
\zeta_{D_n}(s) &=& 1 + n \cdot 2^{-s} + \sum_{d \mid n, \, d \ne 1} (d^{-s} + (2d)^{-s})\\
&=&
(n-1)\cdot 2^{-s} + (1 + 2^{-s}) \prod_{p \mid n}\, f_p(p^{-s}), 
\end{eqnarray*}
where
$$
f_p(X) = X^{m_p} + X^{m_p-1} + \cdots + X + 1, \quad n = \prod_{p\mid n}\, p^{m_p}.
$$ 

\slit
We construct some $p$-groups which are counter examples for Problems [A] and [B].
We denote by $C_n$ the cyclic group of order $n$.
For a prime $p$ and positive integers $m$ and $n$ such that $m > n$, we define a non-abelian group of order $p^{m+n}$ as follows:
\begin{eqnarray}
G_p(m,n) = \gen{\sigma, \tau \mid \sigma^{p^m} = \tau^{p^n} = 1, \; \tau\sigma\tau^{-1} = \sigma^{1+p^{m-1}} }.
\end{eqnarray}

\begin{prop} \label{p-group data}
Assume $G_p(m,n)$ is given as above and $m \geq 3$ if $p = 2$. 

\mmslit
{\rm (1)} The non-abelian group $G_p(m,n)$ and the abelian group $C_{p^m} \times C_{p^n}$ have the same zeta function, and it satisfies the symmetry (0.3). 

\mmslit
{\rm (2)} The $p$-groups of order less than or equal to $p^2$ if $p$ is odd (resp. $2^3$ if $p = 2$) are determined by their zeta functions. 
\end{prop}

\proof
(1) Let us write $G = G_p(m,n)$ and assume $m \geq 3$ if $p = 2$. Since we see
\begin{eqnarray}
(\sigma^i\tau^j)^r = \sigma^{i(r + \frac{r(r-1)}{2}jp^{m-1})}\tau^{rj}, \quad i, j, r \in \N,
\end{eqnarray}
we have
\begin{eqnarray*}
(\sigma^i\tau^j)^p = \left\{\begin{array}{ll}
\sigma^{pi}\tau^{pj} &\mbox{if }\; p \ne 2,\\
\sigma^{2i(1+j2^{m-2})}\tau^{2j} &\mbox{if }\; p = 2.
\end{array} \right.
\end{eqnarray*}
Thus we see 
\begin{eqnarray}
\ord(\sigma^i\tau^j) = p^m \; \mbox{if and only if} \;  p {\not\vert}\, i. 
\end{eqnarray}
If $p\nvert ij$, the group $\gen{\sigma^i\tau^{p^s j}}$ contains an element $\sigma^{i'}\tau^{p^s}$ for some $i'$ not divided by $p$. Further, for $i, \ell$ not divided by $p$, we see that    
$\gen{\sigma^i\tau^{p^s}} = \gen{\sigma^\ell\tau^{p^s}}$ if and only if $i \equiv \ell \pmod{p^{n-s}}$. 
Hence all the subgroups of $G$ isomorphic to $C_{p^m}$ are given as follows:
$$
\gen{\sigma^{i}\tau^{p^s}}, \quad 0 \leq s \leq n, \: 1 \leq i \leq p^{n-s}, \; p\nvert i.
$$ 
For $i$ not divided by $p$, setting $\wt{\sigma} = \sigma^i\tau^{p^s}$, we have $G = \gen{\wt{\sigma}, \tau}$ with the same relations for $\sigma$ and $\tau$. 
There is a natural bijective correspondence between the set of subgroups of $G$ containing $\gen{\wt{\sigma}}$ and that of $G/\gen{\wt{\sigma}} \cong C_{p^n}$, and there is only one subgroup $K_r$ of $G$ such that $(K_r : \gen{\wt{\sigma}}) = p^r$ for each $r$ with $0 \leq r \leq n$.  

On the other hand, any subgroup of $G$ without elements of order $p^m$ is contained in $\gen{\sigma^p, \tau} \cong C_{p^{m-1}} \times C_{p^n}$.      

Hence we see that the subgroups of $G$ naturally and bijectively correspond to those of $C_{p^m} \times C_{p^n}$, hence
their zeta functions coincide.

\mmslit
(2) It is well-known that any $p$-group of order less than or equal to $p^2$ is isomorphic to $C_p, C_{p^2}$ or $C_p \times C_p$, and any $2$-group of order $8$ is isomorphic to one of the following groups
\begin{eqnarray*}
&&
C_8, \quad C_4 \times C_2, \quad C_2 \times C_2 \times C_2, \\
&&
D_4 = \gen{\sigma, \tau \mid \sigma^4 = \tau^2 = 1, \; \tau\sigma\tau = \sigma^{-1}} \left(=G_2(2,1) \right),\\
&&
Q_2 = \gen{\sigma, \tau \mid \sigma^4 = 1, \; \sigma^2 = \tau^2, \; \tau\sigma\tau^{-1} = \sigma^{-1}}.
\end{eqnarray*}
By a direct calculation, we have
\begin{eqnarray} 
&&
\zeta_{C_p}(s) = 1 + p^{-s}, \nonumber \\
&&
\zeta_{C_{p^2}}(s) = 1 + p^{-s} + p^{-2s},\nonumber \\
&& \label{list small odd p}
\zeta_{C_p\times C_p}(s) = 1 + (p+1) p^{-s} + p^{-2s},\nonumber \\
&&
\zeta_{C_8}(s) = 1 + 2^{-s} + 2^{-2s} + 2^{-3s}, \nonumber \\
&&
\zeta_{C_4 \times C_2} = 1 + 3\cdot 2^{-s} + 3\cdot 2^{-2s} + 2^{-3s},\nonumber \\
&&
\zeta_{C_2 \times C_2 \times C_2}(s) = 1 + 7\cdot 2^{-s} + 7\cdot 2^{-2s} + 2^{-3s},\nonumber \\
&&
\zeta_{D_4}(s) = 1 + 5\cdot 2^{-s} + 3\cdot 2^{-2s} + 2^{-3s}, \nonumber \\
&&
\zeta_{Q_2}(s) = 1 + 2^{-s} + 3\cdot 2^{-2s} + 2^{-3s}.  \label{2-groups}
\end{eqnarray}
\qed  

\begin{rem}{\rm
As for the number $a_{p^k}(G_p(m,n)) = a_{p^k}(C_{p^m} \times C_{p^n})$, one can easily obtain directly, or by using Lemma~\ref{lem by Steh} in \S 2 for $\lam = (m,n) \in \Lam_2^+$ and $N_k(\lam) = a_{p^k}(C_{p^m} \times C_{p^n})$.
} %rm
\end{rem}

\slit
For a positive integer $n$ and a prime $p$, we denote by $v_p(n)$ the maximal exponent of $p$ dividing $n$. 
By Propositions~\ref{prop Euler prod} and \ref{p-group data}, we obtain the next theorem, which gives a negative answer for Problems [A] and [B].
We note here that the groups $D_4$ and $Q_2$ are not abelian and their zeta functions are not symmetric (cf. (\ref{2-groups})).
 
\begin{thm}
{\rm (1)} 
For a positive integer $n$ such that $v_p(n) \geq 3$ for some odd prime $p$ or $v_2(n) \geq 4$, there exist abelian groups and non-abelian nilpotent groups of order $n$ which have the same zeta functions. 

\mmslit
{\rm (2)} 
For a positive integer $n$ such that $v_p(n) \leq 2$ for each odd prime $p$ and $v_2(n) \leq 3$, 
any nilpotent group $G$ of order $n$ is determined up to isomorphism by its zeta function $\zeta_G(s)$, and it satisfies the symmetry (\ref{sym of zeta}) if and only if $G$ is abelian. 
\end{thm}

\vspace{2cm}
\Section{}
In this section we enumerate the number of subgroups for abelian $p$-groups. It is well-known that every finite abelian $p$-group is isomorphic to one of $G_\lam$ of the following type:
\begin{eqnarray} 
&&
G_\lam = \Z/p^{\lam_1}\Z \oplus \cdots \oplus \Z/p^{\lam_n}\Z,  \quad \lam \in \Lam_n^+
\end{eqnarray}
where
\begin{eqnarray}
&&
\Lam^+ =   \bigsqcup_{n \geq 1}\, \Lam_n^+,  \qquad \Lam_n^+ = \set{\lam \in \Z^n}{\lam_1 \geq \lam_2 \geq \cdots \geq \lam_n > 0}.
\end{eqnarray}
Hereafter we write $\zeta_\lam(s)$ instead of $\zeta_{G_\lam}(s)$.
We denote by $N_k(\lam)$ the number of subgroups of $G_\lam$ of order $p^k$.
We note $N_k(\lam) = 0$ if $k < 0$ or $k > \abs{\lam}$, for any $\lam \in \Lam^+$.

\slit
The following is clear for $\lam \in \Lam_n^+$:
\begin{eqnarray} 
&&
\abs{G_\lam} = p^{\abs{\lam}}, \quad \abs{\lam} = \sum_{i=1}^n\, \lam_i,\\
&&  
\zeta_\lam(s) = \sum_{k=0}^{\abs{\lam}}\, N_k(\lam) p^{-ks}, \\
&& \label{N-1 rank}
N_1(\lam) = (p^n-1)/(p-1) = p^{n-1} + p^{n-2} + \cdots + p + 1, \\[2mm]
&& \label{symmetry}
N_k(\lam) = N_{\abs{\lam}-k}(\lam), \; 0 \leq k \leq \abs{\lam},
\end{eqnarray}
where (\ref{symmetry}) is a restatement of the duality of abelian groups (cf. (\ref{sym of zeta})). 
Moreover, it is known that 
$N_k(\lam)$ is a polynomial in $p$ with integral coefficients and unimodal, i.e., $N_k(\lam) - N_{k-1}(\lam)$ is a polynomial in $p$ with nonnegative coefficients for $1 \leq k \leq \abs{\lam}/2$ (cf. \cite{But1}, Remark~\ref{rem1}). 

\slit
We prepare some notations.
Let $\lam \in \Lam_n^+$.
Define
\begin{eqnarray}
&&
c_\ell(\lam) = \lam_{\ell+1} + \lam_{\ell+3} + \cdots + 
\left(\begin{array}{ll} \lam_n & \mbox{if}\; n \,{\not\equiv}\, \ell\pmod{2}\\
\lam_{n-1} & \mbox{if}\; n \equiv \ell\pmod{2} \end{array} \right), \; 1 \leq \ell \leq n-1,\nonumber \\[2mm]
&&
c_n(\lam) = 0, %\quad c_{-n}(\lam) = \abs{\lam},  
\quad
c_{-\ell}(\lam) = \abs{\lam} - c_\ell(\lam), 
\end{eqnarray}
in particular
\begin{eqnarray}
ev_\lam = c_1(\lam) = \displaystyle{\sum_{i=1}^{[n/2]}} \lam_{2i}, \quad 
od_\lam =  c_{-1}(\lam) = \displaystyle{\sum_{i=0}^{[n-1/2]}} \lam_{2i+1},
\end{eqnarray}
and we divide the interval  $[0, \abs{\lam}]$ into $(2n-1)$ small intervals $J_\ell(\lam)$ as follows. 
\begin{eqnarray}
&&
J_0(\lam) = [ev_\lam, od_\lam], \\
&&
J_\ell(\lam) =  [c_{\ell+1}(\lam), c_{\ell}(\lam)],\quad
J_{-\ell}(\lam) = [c_{-\ell}(\lam), c_{-(\ell+1)}(\lam)], \quad (1\leq \ell \leq n-1). \nonumber
\end{eqnarray}
We also define
\begin{eqnarray} \label{a-l(lam)}
&&
a_{\ell}(\lam) = \sum_{i=\ell+2}^n \left[\frac{i-\ell}{2}\right] \lam_i, \; (0 \leq \ell \leq n-1), 
\end{eqnarray}
where we understand $a_{n-1}(\lam)  = 0$.

\slit
For $\lam \in \Lam_n^+, \; (n \geq 2)$, we define
\begin{eqnarray}
\lam^{(\ell)} = (\lam_{\ell+1}, \ldots, \lam_{n-1}, \lam_n) \in \Lam_{n-\ell}^+, \quad 1 \leq \ell \leq n-1,
\end{eqnarray}
and we set  $\lam' = \lam^{(1)}$ for simplicity. It is easy to see
\begin{eqnarray}
J_\ell(\lam) = J_0(\lam^{(\ell)}), \quad \ell \geq 1.
\end{eqnarray}

\slit
For a polynomial $g(t)$ in $t$, we denote by $\hterm(g(t))$, $\hdeg(g(t))$ and $\htermcoef(g(t))$ its highest term, highest degree and highest coefficient in $t$, in order. 
We recall $N_k(\lam)$ can be regarded as a polynomial of $p$, and set \\
$H_n(\lam) = \max\set{\hdeg(N_k(\lam))}{0 \leq k \leq \abs{\lam}}$.

%(By modality, $H_n(\lam) = a_0(\lam)$. We calculate directly)  

\begin{thm}\label{thm main key} 
The highest term of $N_k(\lam), \; \lam \in \Lam_n^+$,  as a polynomial in $p$ is given as
\begin{eqnarray*}
&&
\hterm(N_k(\lam)) = C(n, \lam; k)p^{\ell k + a_\ell(\lam)}, \quad \mbox{if }\; k \in J_\ell(\lam), \; 0 \leq \ell \leq n-1,\\[2mm]
&&
\hterm(N_k(\lam)) = \hterm(N_{\abs{\lam}-k}(\lam)), \quad \mbox{if }\, \abs{\lam}/2 \leq k \leq \abs{\lam}\\
&&
\hspace*{2.4cm}  \left(= C(n, \lam; \abs{\lam}-k)p^{-\ell k + \ell \abs{\lam}+a_\ell(\lam)}, \quad \mbox{if }k \in J_{-\ell}(\lam), \; \ell \geq 1\right).
\end{eqnarray*}
Here $C(n, \lam; k) = 1$ if $n \leq 2$, and in general 
\begin{eqnarray}
\lefteqn{C(n, \lam; k)} \nonumber \\
&=&
\left\{\begin{array}{ll}
C(n-\ell, \lam^{(\ell)}; k) \qquad & \mbox{if }\; k \in J_\ell(\lam) = J_0(\lam^{(\ell)}), \; \ell \geq 1,\\[2mm]
\displaystyle{ \sum_{i=\max\{ev_\lam, \, k-\lam_1+\lam_2\}}^{\min\{k, \, od_\lam-\lam_1+\lam_2\}} }C(n-1,\lam'; i)
& \mbox{if }\; k \in J_0(\lam) .
\end{array}\right.
\end{eqnarray}
Further, for each $\ell \geq 0$, $C(n, \lam; k), k \in J_\ell(\lam)$ is unimodal and symmetric with respect to ${\abs{\lam^{(\ell)}}}/{2}$ and 
$C(n, \lam; c_\ell(\lam)) = 1$.
In particular, $H_n(\lam) = a_0(\lam) = \sum_{i=2}^n \left[\frac{i}{2}\right] \lam_i$.
\end{thm}

\mslit
\begin{rem}{\rm 
We see the leading coefficients of smaller sizes appear repeatedly like a fractale, since $C(n, \lam; k) = C(n-\ell, \lam^{(\ell)}; k)$ for $k \in J_\ell(\lam)$ with $\ell \geq 1$. The new coefficients appear for $k \in J_0(\lam)$.

\mmslit
It may happen $J_\ell(\lam) = \{c_\ell(\lam)\}$, but $c_\ell(\lam) > c_{\ell+2}(\lam)$ for any $\ell$, $-n \leq \ell \leq n-2$. 

\mmslit
We give some examples of $C(n, \lam; k)$ for $k \in J_0(\lam)$:   
\begin{eqnarray}
C(3, \lam ;k) 
&=& \left\{\begin{array}{ll}
k-\lam_2+1 & \mbox{if }\; k \in [\lam_2, \lam_2+\lam_3]\\
\lam_3+1 & \mbox{if }\; k \in [\lam_2+\lam_3, \lam_1]\\
-k+\lam_1+\lam_3+1 & \mbox{if }\; k \in [\lam_1, \lam_1+\lam_3] ,
\end{array} \right. \nonumber \\
\mbox{or} && \\
&=& \left\{\begin{array}{ll}
k-\lam_2+1 & \mbox{if }\; k \in [\lam_2, \lam_1]\\
\lam_1-\lam_2+1 & \mbox{if }\; k \in [\lam_1, \lam_2+\lam_3]\\
-k+\lam_1+\lam_3+1 & \mbox{if }\; k \in [\lam_2+\lam_3,  \lam_1+\lam_3] .
\end{array} \right. \nonumber \\[2mm]
C(4, \lam ;k) 
&=& \left\{\begin{array}{ll}
k-\lam_2-\lam_4+1 & \mbox{if }\; k \in [\lam_2+\lam_4, \lam_2+\lam_3]\\
\lam_3-\lam_4+1 & \mbox{if }\; k \in [\lam_2+\lam_3, \lam_1+\lam_4]\\
-k+\lam_1+\lam_3+1 & \mbox{if }\; k \in [\lam_1+\lam_4, \lam_1+\lam_3] ,
\end{array} \right. \nonumber \\
\mbox{or} && \\
&=& \left\{\begin{array}{ll}
k-\lam_2-\lam_4+1 & \mbox{if }\; k \in [\lam_2+\lam_4, \lam_1+\lam_4]\\
\lam_1-\lam_2+1 & \mbox{if }\; k \in [\lam_1+\lam_4, \lam_2+\lam_3]\\
-k+\lam_1+\lam_3+1 & \mbox{if }\; k \in [\lam_2+\lam_3,  \lam_1+\lam_3] .
\end{array} \right. \nonumber 
\end{eqnarray}
} %rm
\end{rem}

\mslit
The following lemma follows from a result of Stehling(\cite{Steh}) and plays a key role for the proof of Theorem~\ref{thm main key}.
\begin{lem} \label{lem by Steh}
For any $\lam \in \Lam_n^+$, set $\lam' = (\lam_2, \ldots, \lam_n) \in \Lam_{n-1}^+$. Then one has
\begin{eqnarray} \label{eq in Lem1}
N_k(\lam) = \sum_{i=0}^k p^i N_i(\lam') -\sum_{i=\abs{\lam}+1-k}^{\abs{\lam'}}\, p^i N_i(\lam'), \quad 0 \leq k \leq \abs{\lam},
\end{eqnarray}
where the second summation appears only when  $k > \lam_1$.
\end{lem}
\proof
For $\mu \in \Lam_n^+$ satisfying $\mu_1= \cdots = \mu_i > \mu_{i+1}$ for some $i \geq 1$, \cite[Theorem 1]{Steh} claims
\begin{eqnarray} \label{Steh Th1}
N_k(\mu) = N_k(\mu^*) + p^{\abs{\mu}-k} N_{k-\mu_i}(\mu^\vee), %\quad (k \geq \mu_i),
\end{eqnarray}
where $\mu^* \in \Lam_n^+$ defined by $\mu^*_i = \mu_i-1$ and $\mu^*_j = \mu_j$ except $j =i$, and \\
$\mu^\vee = (\mu_1, \ldots, \mu_{i-1}, \mu_{i+1}, \ldots, \mu_n) \in \Lam_{n-1}^+$. 
For any $\lam \in \Lam_n^+$, set $\wt{\lam} = (\lam_1+1, \lam_2, \ldots, \lam_n) \in \Lam_n^+$ and $\lam' = (\lam_2, \ldots, \lam_n) \in \Lam_{n-1}^+$. Then, since ${\wt{\lam}}^* = \lam$ and ${\wt{\lam}}^\vee = \lam'$, one has by (\ref{Steh Th1})
\begin{eqnarray} \label{(1)}
N_j(\wt{\lam}) = N_j(\lam) + p^{\abs{\lam}+1-j} N_{j-\lam_1-1}(\lam'). 
\end{eqnarray}
By the duality (\ref{symmetry}), and again by (\ref{Steh Th1}),
\begin{eqnarray}
N_j(\wt{\lam}) &=&
N_{\abs{\lam}+1-j}(\wt{\lam}) = N_{\abs{\lam}+1-j}(\lam) + p^j N_{\abs{\lam'}-j}(\lam') \nonumber\\
&=& \label{(2)}
N_{j-1}(\lam) + p^j N_j(\lam').
\end{eqnarray}
Summing up RHD's of (\ref{(1)}) and (\ref{(2)})  from $j = 1$ to $k$, one has (the terms $N_j(k)$ is cancelled for $1 \leq j \leq k-1$)
\begin{eqnarray*}
N_k(\lam) + \sum_{j = 1}^k\, p^{\abs{\lam}+1-j} N_{j-\lam_1-1}(\lam') = 
N_0(\lam) + \sum_{j = 1}^k\, p^j N_j(\lam'),
\end{eqnarray*}
hence (setting $i = \abs{\lam}+1-j$ in the LHS of the above), we have
\begin{eqnarray*}
N_k(\lam) &=&
\sum_{j = 0}^k\, p^j N_j(\lam') -  \sum_{i = \abs{\lam}+1-k}^{\abs{\lam}}\, p^i N_{\abs{\lam}-\lam_1-i}(\lam') \\
&=&
\sum_{j = 0}^k\, p^j N_j(\lam') -  \sum_{i = \abs{\lam}+1-k}^{\abs{\lam'}}\, p^i N_{\abs{\lam'}-i}(\lam')\\
&=&
\sum_{j = 0}^k\, p^j N_j(\lam') -  \sum_{i = \abs{\lam}+1-k}^{\abs{\lam'}}\, p^i N_i(\lam'),\\
\end{eqnarray*}
where the second summation appears only when $\abs{\lam}+1-k \leq \abs{\lam'}$, i.e. $k > \lam_1$. 
\qed

\slit
For $r \in \Lam_1^+$, since $G_r$ is a cyclic group of order $p^r$, one has
\begin{eqnarray} \label{n=1}
N_k(r) = 1, \quad \mbox{for }\; k \in J_0(r) = [0, r]. 
\end{eqnarray}
Then, for $\mu \in \Lam_2^+$, one has by Lemma~\ref{lem by Steh}
\begin{eqnarray} \label{n=2}
\hterm(N_k(\mu)) = \left\{ \begin{array}{ll}
p^k & \mbox{for }\; k \in J_1(\mu) = [0, \mu_2]\\
p^{\mu_2} & \mbox{for }\; k \in J_0(\mu) = [ev_\mu, od_\mu] = [\mu_2, \mu_1]\\
p^{\abs{\mu}-k} & \mbox{for }\; k \in J_{-1}(\mu) = [\mu_1, \abs{\mu}].
\end{array}
\right. 
\end{eqnarray}

\slit
We note the following useful relations :
\begin{eqnarray} 
&&
\label{rel-1}
c_{-1}(\lam') = \abs{\lam'}-c_1(\lam') = ev_\lam, \quad c_{\ell-1}(\lam') = c_\ell(\lam),\; (\ell \geq 1)\\
&& 
\label{rel-a}
a_{\ell-1}(\lam') = a_\ell(\lam), \quad (\ell \geq 1),\\
&&\label{a-0(lam)}
a_0(\lam) = \lam_2 + \lam_3 + 2\lam_4+2\lam_5 + 3\lam_6 + 3\lam_7 + \cdots = \abs{\lam'} + a_1(\lam')\\
&& \label{rel a to c}
a_\ell(\lam) + c_\ell(\lam) = a_{\ell-1}(\lam).
\end{eqnarray}

\slit

\begin{lem}  \label{lem: induction}
Assume $n \geq 3$ and Theorem~\ref{thm main key} is established for $\Lam_{n-1}^+$, and take any $\lam \in \Lam_n^+$.
Then the highest degree $\hdeg(p^iN_i(\lam'))$ is strictly monotone increasing for $i \in [0, ev_\lam]$, 
stable for $i \in [ev_\lam, od_\lam-\lam_1+\lam_2]$ with value $a_0(\lam)$, and strictly monotone decreasing for $i \in [od_\lam-\lam_1+\lam_2, \abs{\lam'}]$. 
\end{lem}

\proof  
%Let $0 \leq \ell \leq n-1$. 
We may use the result of Theorem~\ref{thm main key} for $\lam' \in \Lam_{n-1}^+$.

\mmslit
For $i \in [0, c_{-1}(\lam')]$, the highest degree $\hdeg(p^iN_i(\lam'))$ is strictly monotone increasing 
 and takes the highest value at $c_{-1}(\lam') = ev_\lam$ with value 
\begin{eqnarray}
c_{-1}(\lam') + a_0(\lam') = c_1(\lam) + a_1(\lam) = a_0(\lam)
\end{eqnarray}
by (\ref{rel-a}) and (\ref{rel a to c}). 

\mmslit
For $i \in J_{-1}(\lam') = [c_{-1}(\lam'), c_{-2}(\lam')]$, the highest degree $\hdeg(p^iN_i(\lam'))$ is stable with value
\begin{eqnarray}
&&
\hdeg(p^iN_i(\lam')) = \abs{\lam'} + a_1(\lam') = a_0(\lam),
\end{eqnarray}
by (\ref{a-0(lam)}), and
\begin{eqnarray}
&&
c_{-2}(\lam') = \abs{\lam'} - c_2(\lam') = od_\lam -\lam_1 + \lam_2.
\end{eqnarray}

\mmslit
For $i \in [c_{-2}(\lam'), \abs{\lam'}]$, the highest degree $\hdeg(p^iN_i(\lam'))$ is strictly monotone decreasing and takes the highest value at $i= c_{-2}(\lam')$ with value
\begin{eqnarray*}
-c_{-2}(\lam')+2\abs{\lam'} + a_2(\lam') = \abs{\lam'} + c_2(\lam') + a_2(\lam') = a_0(\lam)
\end{eqnarray*}
by (\ref{a-0(lam)}) and (\ref{rel a to c}).
\qed

\vspace{1cm}
%Assume $0 \leq k \leq \abs{\lam'}$.
%
\newcommand{\negativeterm}{{\sum_{i=\abs{\lam}+1-k}^{\abs{\lam'}}\, p^i N_i(\lam')}}
\slit
{\it Proof of Theorem~\ref{thm main key}}. \;  
We assume that the statement holds for $n-1$, and take any $\lam \in \Lam_n^+$ and fix it.
We set 
\begin{eqnarray}
(pos)_k = \sum_{i=0}^k\, p^i N_i(\lam'), \quad  (neg)_k = \negativeterm
\end{eqnarray}
where $(neg)_k$ appears only when $k > \lam_1$.
By Lemma~\ref{lem: induction}, we see 
\begin{eqnarray}
\label{pos-hi}
\hdeg((pos)_k) = a_0(\lam) &\Longleftrightarrow& k \geq ev_\lam,\\
\hdeg((neg)_k) = a_0(\lam) &\Llra& k > \lam_1 \, \mbox{and}\, \abs{\lam}-k+1 \leq od_\lam-\lam_1+\lam_2
 \nonumber\\
{}& \label{neg-hi}
\Llra& k \geq ev_\lam + \lam_1-\lam_2+1.
\end{eqnarray}
Here we note that the condition (\ref{neg-hi}) yields $k > \lam_1$ and $k > ev_\lam$.

%If $k \leq \lam_1$, it is clear that $\hterm(N_k(\lam)) = \hterm((pos)_k)$.

\mslit
Case 1 [$k \leq ev_\lam$] \\ 
There exists an $\ell \geq 1$ such that $k \in J_{\ell-1}(\lam') = J_\ell(\lam)$ and  
\begin{eqnarray}
\hterm((pos)_k) 
&=&  
\hterm(p^kN_k(\lam')) = C(n-1, \lam'; k)p^{\ell k + a_{\ell-1}(\lam')}  \nonumber\\
&=& \label{posdeg case1}
C(n-1, \lam'; k)p^{\ell k + a_{\ell}(\lam)}.
\end{eqnarray}
If $k \leq \lam_1$, (\ref{posdeg case1}) coincides with $\hterm(N_k(\lam))$. Assume $k > \lam_1$. Then, since $c_{\ell-1}(\lam') \geq k$, 
\begin{eqnarray*}
\abs{\lam}-k + 1 \geq \abs{\lam}-c_{\ell-1}(\lam')+ 1 = \abs{\lam'}-(c_{(\ell-1)}(\lam')-\lam_1 - 1) \geq  c_{-(\ell+1)}(\lam'),
\end{eqnarray*}
%($+(\lam_1+1)$ causes the move of the interval at least $+2$)
and
\begin{eqnarray}
\lefteqn{\hdeg((neg)_k)) = \hdeg(p^{\abs{\lam}-k+1}N_{\abs{\lam}-k+1}(\lam')) } \label{negdeg case1}\\
&\leq& -\ell(\abs{\lam}-k+1) + (\ell+1)\abs{\lam'} + a_{\ell+1}(\lam') %\nonumber \\
=  
-\ell \lam_1+\abs{\lam'}+\ell(k-1)+a_{\ell+2}(\lam). \nonumber
\end{eqnarray}
Hence, by (\ref{posdeg case1}) and (\ref{negdeg case1}), we have 
\begin{eqnarray*}
\hdeg((pos)_k)) - \hdeg((neg)_k)) &\geq& a_{\ell}(\lam) - a_{\ell+2}(\lam) - \abs{\lam'} + \ell \lam_1 + \ell \nonumber\\
&=&
-(\lam_2+\cdots + \lam_{\ell+1})+\ell \lam_1 + \ell > 0.  
\end{eqnarray*} 
Hence,  independently whether $k > \lam_1$ or not, we have for $k \leq ev_\lam$,
\begin{eqnarray}  \label{Case1 result}
\hterm(N_k(\lam)) &=& \hterm(p^kN_k(\lam')) \nonumber\\
&=&
C(n-1, \lam'; k)p^{\ell k + a_{\ell}(\lam)} \quad (k \in J_\ell(\lam), \ell \geq 1). 
\end{eqnarray}
Further, if $k \in J_\ell(\lam)$ with $\ell \geq 2$, then we see
\begin{eqnarray}
C(n, \lam; k)  &=& C(n-1, \lam'; k) = C(n-1-(\ell-1), \lam'^{(\ell-1)}; k) \nonumber \\
& = &
C(n-\ell, \lam^{(\ell)}; k),
\end{eqnarray}
and in particular, 
\begin{eqnarray}  \label{edge=1}
C(n, \lam; c_\ell(\lam)) &=&
C(n-1, \lam'; c_{\ell-1}(\lam')) = 1, \quad \ell \geq 2.
\end{eqnarray}

\mslit
Case 2 [$k \geq od_\lam$]\\
Since $N_k(\lam) = N_{\abs{\lam}-k}(\lam)$ and $\abs{\lam}-k \leq ev_\lam$, the result follows from Case1.
Thus, for some $\ell \geq 1$, $k \in J_{-\ell}(\lam)$ and $\abs{\lam}-k \in J_\ell(\lam)$, and 
\begin{eqnarray}
\hterm(N_k(\lam)) &=& 
C(n, \lam; \abs{\lam}-k) p^{\ell(\abs{\lam}-k)+a_\ell(\lam)} \nonumber\\
&=&
C(n-\ell, \lam^{(\ell)}; \abs{\lam}-k) p^{\ell(\abs{\lam}-k)+a_\ell(\lam)}.
\end{eqnarray}

\mslit
Case 3 [$ev_\lam \leq k \leq ev_\lam+\lam_1-\lam_2$] \\
By (\ref{pos-hi}) and (\ref{neg-hi}), we have
\begin{eqnarray} \label{Case3 result}
&&
\hdeg(N_k(\lam)) = \hdeg((pos)_k) = a_0(\lam), \nonumber \\[2mm]
&& \label{case3-coeff}
\htermcoef(N_k(\lam)) = \sum_{i=ev_\lam}^{\min\{k, \, od_\lam-\lam_1+\lam_2 \}} C(n-1,\lam'; i). 
\end{eqnarray}

\slit
Case 4 [$ev_\lam+\lam_1-\lam_2+1 \leq k \leq od_\lam$] \\
We see $\hdeg((pos)_k) = \hdeg((neg)_k) = a_0(\lam)$ and $\abs{\lam}-k+1 \geq ev_\lam+1$, and the difference of the highest coefficients becomes
\begin{eqnarray} \lefteqn{\htermcoef((pos)_k)  - \htermcoef((neg)_k)} \nonumber\\
&=& \label{Case4-1}
\sum_{i=ev_\lam}^{\min\{k, \, od_\lam-\lam_1+\lam_2\}} C(n-1,\lam'; i)
-
\sum_{i = \abs{\lam}-k+1}^{od_\lam-\lam_1+\lam_2} C(n-1,\lam'; i).
\end{eqnarray}
Since $C(n-1,\lam'; i)$ is symmetric for $i \in J_{-1}(\lam') = [ev_\lam, od_\lam-\lam_1+\lam_2]$ with respect to $(\abs{\lam}-\lam_1+\lam_2)/2$ and 
\begin{eqnarray*}
od_\lam-\lam_1+\lam_2-\abs{\lam}+k &=& k-ev_\lam-\lam_1+\lam_2 \\
&<&
\min\{k, \, od_\lam-\lam_1+\lam_2\}-ev_\lam+1
\end{eqnarray*}
we see (\ref{Case4-1}) becomes 
\begin{eqnarray} 
\label{Case4-2}
\sum_{i=k-\lam_1+\lam_2}^{\min\{k, \, od_\lam-\lam_1+\lam_2\}} C(n-1,\lam'; i) > 0,
\end{eqnarray}
and $\hdeg(N_k(\lam)) = a_0(\lam)$ with coefficient (\ref{Case4-2}).

\mslit
By Case 3 and Case 4, we obtain for $k \in J_0(\lam) = [ev_\lam, od_\lam]$, 
\begin{eqnarray}
&&
\hterm(N_k(\lam)) = C(n,\lam;k) p^{a_0(\lam)}, \nonumber\\[2mm]
&&   \label{C-n-lam}
C(n,\lam;k) = \sum_{\max\{ev_\lam, \, k-\lam_1+\lam_2\}}^{\min\{k, \, od_\lam-\lam_1+\lam_2\}}\, C(n-1,\lam'; k).
\end{eqnarray}
By (\ref{C-n-lam}) and the induction hypothesis, the unimodality of $C(n, \lam; k)$ for $k \in J_0(\lam)$ is clear, and one has for $c_1(\lam) = ev_\lam = c_{-1}(\lam')$, 
\begin{eqnarray}
C(n, \lam; c_1(\lam)) = C(n-1, \lam'; c_1(\lam)) = C(n-1, \lam', c_{-1}(\lam')) = 1,
\end{eqnarray}
which completes the proof of Theorem~\ref{thm main key}.
\qed

\slit
As a corollary of Theorem~\ref{thm main key}, we obtain a result for zeta functions, where we don't need the information of the leading coefficients, that of the highest degrees is enough.

\begin{thm} \label{thm lam=mu}
The zeta function $\zeta_\lam(s)$ determines $\lam \in \Lam^+$.
\end{thm}

\proof
The zeta function $\zeta_\lam(s)$ determines $\abs{\lam}$ and $n$ for which $\lam \in \Lam_n^+$ (cf. (\ref{N-1 rank})).  
Take $\lam, \mu \in \Lam_n^+$ such that $\zeta_\lam(s) = \zeta_\mu(s)$, then $\hdeg(N_k(\lam)) = \hdeg(N_k(\mu))$ at any $k$. Since the statement for $n \leq 2$ is clear by (\ref{n=1}) and (\ref{n=2}), we consider the case $n \geq 3$.  
Assume $\lam_n < \mu_n$. Taking $k = \lam_n+1$, we have $(n-2)(\lam_n+1)+\lam_n = (n-1)(\lam_n+1)$ by Theorem~\ref{thm main key}, which is a contradiction; thus $\lam_n = \mu_n$. In the same way, we obtain $\lam_i = \mu_i$ from $i=n-1$ to $i =1$ in order. Thus $\lam = \mu$.
\qed

\slit
Since any finite abelian group is decomposed into the direct product of its Sylow $p$-subgroups, combining Theorem~\ref{thm lam=mu} and Proposition~1.2, we obtain the affirmative statement for Problem [A]: 

\begin{thm} \label{thm answer for [A]}
Within finite abelian groups, the zeta function $\zeta_G(s)$ determines the isomorphism class of $G$.
\end{thm}

\slit
\begin{rem}\label{rem1} {\rm 
The enumerating subgroups of finite abelian groups is a classical problem(cf. \cite{But1}, \cite{But2}). The number $g^{\lam}_{\mu \nu}(p)$ of subgroups $H$ of $G_\lam$ such that $H \cong G_\nu$ and $G_\lam/H \cong G_\mu$ was studied by G.~Hall in 1950's, and it is known to have integral coefficients as a polynomial in $p$. These numbers $g^{\lam}_{\mu \nu}(p)$ appear as coefficients of the product of Hall-Littlewood symmetric polynomials $P_\lam(x; t)$, which can be considered as a system of orthogonal polynomials associated with the root system of type $A_n$, and I.~G.~Macdonald have introduced and studied more general orthogonal polynomials associated with various root systems(cf. \cite{Mac}, \cite{Mac2}). The unimodality of $N_k(\lam)$ is proved by L.~M.~Butler by using a result on Hall-Littlewood polynomials(\cite{But1}). 
}%rm
\end{rem}

\begin{rem}\label{rem2}{\rm 
The zeta functions $\zeta_\lam(s)$ has a close relation to Eisenstein series of $GL_{2n}$ and local densities of square matrices, which we explain after \cite{FS}. Let $T \in M_n(\Z)$ of rank $n$ and a prime $p$. Define
\begin{eqnarray*}
b_p(T, s) = \int_{M_n(\Q_p)} \nu_p(X)^{-s}{\bf e}_p(2\pi i \, \tr({}^tTX))dX,
\end{eqnarray*} 
where $dX$ is the Haar measure on $M_n(\Q_p)$, $\nu_p(X)$ is the product of the $p$-denominators of the elementary divisors of $X$ and ${\bf e}_p(x) = \exp(2\pi i x')$, $x'$ being the fractional part. Then $b(T, s) = \prod_{p} b_p(T, s)$ appears as a main term in the Fourier expansion of $E_{n,n}(Y, s)$ for symmetric $Y \in M_{2n}(\Z)$, and $b(T, s)$ is an analogue of the Siegel singular series.  
The local density $\alp_p(I_m, T)$ is defined by 
\begin{eqnarray*}
\alp_p(I_m, T) = \lim_{\ell \rightarrow \infty}\, \frac{1}{ p^{\ell(m^2-n^2)}} \sharp\set{(X, Y) \in M_{mn}(\Z/p^\ell \Z)^{\oplus 2}}{{}^tXY \equiv T \mod {p^\ell}}.
\end{eqnarray*}
$b_p(T;s)$ and $\alp_p(I_m, T)$ are determined by the two-sided $GL_n(\Z_p)$-coset containing $T$, which is represented by some $D_\lam = Diag(p^{\lam_1}, \ldots, p^{\lam_n})$ with $\lam \in \wt{\Lam_n^+} = \set{\lam \in \Z^n}{\lam_1 \geq \cdots \geq \lam_n \geq 0}$. Then it is known (cf. \cite[Theorem~2]{FS}, also \cite{BB}) 
\begin{eqnarray*}
&& 
b_p(D_\lam,  s) = \zeta_\lam(s-n) \cdot \prod_{i=0}^{n-1}\, (1-p^{-(s-i)}), \\
&&
b_p(D_\lam, m) = \alp_p(I_m, D_\lam), \quad (m \geq n \geq 1). 
\end{eqnarray*}  
Here $\lam$ can be regarded in $\Lam^+ = \cup_n \Lam_n^+$ unless $\lam = (0, \ldots, 0)$, which corresponds to the trivial group $\{e\}$ and $\zeta_{\{e\}}(s) = 1$. By Theorem~\ref{thm lam=mu}, we see $b_p(D_\lam, s)$ is different for each $\lam \in \wt{\Lam_n^+}$. 
} % rm
\end{rem}

\vspace{2cm}

\Section{}
We keep the notation in \S 2, especially $G_\lam$ indicates the abelian $p$-group of order $p^{\abs{\lam}}$ defined in (2.1) for each $\lam \in \Lam^+$. %(\ref{}).
We consider whether there is a group that is not isomorphic to $G_\lam$ but has the same zeta function with $G_\lam$ or not.

\begin{prop}
{\rm (1)} Let $\lam \in \Lam_n^+$ and assume $n = 1$ or $\lam_1 =1$. For a group $G$, the coincidence $\zeta_G(s) = \zeta_\lam(s)$ implies the isomorphism $G \cong G_\lam$.

{\rm (2)} Let $\lam \in \Lam_n^+$ with $n \geq 2$. Assume  $\lam_1 > \lam_2$ and $\lam_1 \geq 3$ if $p =2$. Then there exists a non-abelian group $\wt{G_\lam}$ for which $\zeta_{\wt{G_\lam}}(s) = \zeta_\lam(s)$, e.g., set
\begin{eqnarray}
\wt{G_\lam} = \langle
a_1, \ldots, a_n  \mid \begin{array}{l}
a_i^{p^{\lam_i}} = 1 \; (1 \leq i \leq n), \; a_na_1a_n^{-1} = a_1^{1+p^{\lam_1-1}}
\\
a_ia_j = a_ja_i \, \mbox{unless}\; \{i, j\} = \{1, n\}
\end{array} \rangle .
\end{eqnarray}
\end{prop}

\proof
(1) The result for the case $n = 1$ follows from Proposition~1.1.
Next, let $\lam = (1^n), \;  n \geq 2$ and assume the zeta function of a group $G$ coincides with $\zeta_\lam(s)$. 
Since 
$$
a_p(G) = N_1(\lam) = \frac{p^n-1}{p-1},
$$
we see the exponent of $G$ is $p$. Denote by $Z(G)$ the center of $G$ and by $Z_G(a)$ the centralizer of $a (\in G)$ in $G$. Then any subgroup of $G$ of order $p^2$ is isomorphic to $C_p \times C_p$, and 
\begin{eqnarray*}
a_{p^2}(G) &\leq& \frac{(\abs{Z(G)}-1)(p^n-p)}{(p^2-1)(p^2-p)} + \sum_{a \in G \backslash Z(G)}\, \frac{\abs{Z_G(a)}-p}{(p^2-1)(p^2-p)}\\
& \leq &
\frac{(p^r-1)(p^n-p) + (p^n-p^r)(p^{n-1}-p)}{(p^2-1)(p^2-p)},
\end{eqnarray*}
where $p^r = \abs{Z(G)}$. Hence
\begin{eqnarray*}
N_2(\lam) - a_{p^2}(G) &\geq&  \frac{(p^n-1)(p^n-p) - \{(p^r-1)(p^n-p) + (p^n-p^r)(p^{n-1}-p)\}}{(p^2-1)(p^2-p)}\\
&=& 
\frac{(p^{2n}-p^{n+r})(1-p^{-1})}{(p^2-1)(p^2-p)}.
\end{eqnarray*}
Thus we see the identity $a_{p^2}(G) = N_2(\lam)$ implies $r = n$, i.e., $G$ is abelian, hence $G \cong G_\lam$ (the elementary abelian group of rank $n$).  

\mmslit
(2) The strategy is similar to the case $G_p(m,n)$ in \S 1,  indeed $G_p(m,n) \cong \wt{G_\lam}$ with $\lam = (m,n) \in \Lam_2^+$. Take any $\lam \in \Lam_n^+$ satisfying the condition of (2) and construct $\wt{G_\lam}$. Any element $x \in \wt{G_\lam}$ can be expressed uniquely as
\begin{eqnarray*}
x = a_1^{e_1}\cdots a_n^{e_n}, \quad 0 \leq e_i < p^{\lam_i} \; (1 \leq i \leq n),
\end{eqnarray*}
and 
\begin{eqnarray*}
ord(x) = p^{\lam_1} \quad \mbox{if and only if} \quad p {\nvert} e_1.
\end{eqnarray*}
Further we see
\begin{eqnarray*}
&&
x^r = a_1^{e_1(r+\frac{r(r-1)}{2}e_np^{\lam_1-1})} \prod_{i \geq 2} a_i^{r e_i}, \quad (r \geq 1), \\
&&
a_nxa_n^{-1} = a_1^{e_1(1+p^{\lam_1-1})} a_2^{e_2} \cdots a_n^{e_n} = x^{1+p^{\lam_1-1}}.
\end{eqnarray*}
Thus, for any $x \in \wt{G_\lam}$ of order $p^{\lam_1}$, $\gen{x}$ is normal in $\wt{G_\lam}$,  
$\wt{G_\lam}/\gen{x} \cong C_{p^{\lam_2}} \times \cdots \times C_{p^{\lam_n}}$, and the set of subgroups of $\wt{G_\lam}$ containing $\gen{x}$ is one to one corresponding to the set of subgroups of $C_{p^{\lam_2}} \times \cdots \times C_{p^{\lam_n}}$.  
On the other hand, any subgroup of $\wt{G_\lam}$ without an element of order $p^{\lam_1}$ is contained in $\gen{{a_1}^p, a_2, \cdots, a_n} \cong C_{p^{\lam_1-1}} \times C_{p^{\lam_2}} \times \cdots \times C_{p^{\lam_n}}$.    
The above situation is similar to the case of $G_\lam$, and we see there is a natural bijection between the set of subgroups of $\wt{G_\lam}$ and that of $G_\lam$, and their zeta functions coincide. 
\qed

\begin{rem}{\rm 
The zeta function of the above group $\wt{G_\lam}$ satisfies the symmetry (0.3), by virtue of its coincidence with that of $G_\lam$. Hence, one has many examples of non-abelian groups within nilpotent groups having zeta functions with symmetry.  
By neumerical data,  we know some examples of non-abelian groups whose zeta functions do not coincide with that of any abelian group but satisfies the symmetry, e.g., the group of GAP ID $(2^4, 10)$ (the smallest order) or GAP ID $(3^4, 14)$ (the smallest odd order).

\mmslit
Let $\lam \in \Lam_n^+$ with $n \geq 2$. Assume $\lam \ne (1, \ldots, 1)$ and also $\lam \ne (2,1\ldots,1)$ if $p = 2$. Then,  we expect there exists a non-abelian $p$-group for which $\zeta_G(s) = \zeta_\lam(s)$. 
}%\rm
\end{rem}

\slit

\slit
Acknowledgement: The author expresses her thanks to A.~Hoshi who checked the groups of order $27$ by using GAP. 
She is thankful also to the members of WINJ, since WINJ7 was a good opportunity to think about this theme.

\vspace{2cm}

\bibliographystyle{amsalpha}

\end{document}